\documentclass[Accepted,10pt]{elsarticle}
\usepackage{amsfonts}
\usepackage{amsmath,amssymb,amscd,enumerate,latexsym,graphicx}

 \usepackage{graphics}
\usepackage{graphicx}
 \usepackage{epsfig}
\usepackage{amssymb}
\usepackage{amsthm}

\newtheorem{thm}{Theorem}

\newdefinition{rmk}{Remark}
\newproof{pf}{Proof}
\newdefinition{example}{Example}
\newdefinition{definition}{Definition}
\newdefinition{property}{Property}
\newdefinition{problem}{Problem}
\newdefinition{corollary}{Corollary}

\journal{L.A.A.}

\begin{document}

\begin{frontmatter}



\title{Exclusion sets for eigenvalues of matrices}


\author[label2]{Suhua Li}
\author[label2]{Chaoqian Li\corref{cor1}}
\ead{lichaoqian@ynu.edu.cn}
\author[label2]{Yaotang Li}
 \cortext[cor1]{Corresponding author.}
\address[label2]{School of Mathematics and Statistics, Yunnan
University, Kunming, P. R. China 650091}

\begin{abstract}
To locate all eigenvalues of a matrix  more precisely, we exclude some sets which
do not include any eigenvalue of the matrix  from the well-known Brauer set to give two new Brauer-type eigenvalue inclusion sets. And it is also shown that  
the new sets are contained in the Brauer set.
\end{abstract}

\end{frontmatter}


\section{Introduction}
One of the most important problems for eigenvalues of matrices is to locate them \cite{H},
that is, to find regions including all eigenvalue of a given matrix $A$ in the complex plane. The well-known Ger\v{s}gorin disk theorem \cite{Ge} stated below provides just such a region,
which consists of $n$ disks centered at the diagonal elements of the matrix.
\begin{thm}\label{th1.1} \emph{\cite{Ge}}
Let $A=[a_{ij}]\in \mathbb{C}^{n\times n}$ be a complex matrix, and $\sigma(A)$ the set of all eigenvalues of $A$. Then
 \[\sigma(A) \subseteq  \Gamma(A) = \bigcup \limits_{i=1}^n  \Gamma_i(A),  \]
 where \[\Gamma_i(A)=\left \{z\in \mathbb{C}:|z-a_{ii}|\leq  r_i(A) \right \},\]
and $r_i(A)= \sum\limits_{k=1,\atop k\neq i}^{n} |a_{ik}|.$
\end{thm}

Generally, $\Gamma_i(A)$ is called the $i$th disk, and $ \Gamma(A)$ is called the Ger\v{s}gorin set.
Although the Ger\v{s}gorin set is beautiful and simple \cite{Ge}, it is only a raw result,
which inspires researchers to find another sets which are tighter than $\Gamma(A)$.
One such well-known set provided by Brauer \cite{Br} is described as follows.

\begin{thm}\label{th1.2} \emph{\cite{Br}}
Let $A=[a_{ij}]\in \mathbb{C}^{n\times n}$ be a complex matrix. Then
 \[\sigma(A) \subseteq  \mathcal{K}(A) = \bigcup \limits_{j\neq i, \atop i,j=1}^{n}  \mathcal{K}_{ij}(A),  \]
where \[\mathcal{K}_{ij}(A)=\left \{z\in \mathbb{C}:|z-a_{ii}||z-a_{jj}|\leq  r_i(A)r_j(A) \right \}.\]
\end{thm}

Note that the Brauer set $\mathcal{K}(A)$ consists of $ \frac{n(n-1)}{2}$ Cassini ovals $\mathcal{K}_{ij}(A)$ \cite{Br}. Hence, $\mathcal{K}(A)$ needs more computations than
 $\Gamma(A)$ to locate all eigenvalues of $A$, while $\mathcal{K}(A)$ can captures all eigenvalues of $A$ more precisely than $\Gamma(A)$, that is
 \[\mathcal{K}(A) \subseteq \Gamma(A).\]

Besides the Brauer set, there are many sets which are all tighter than the Ger\v{s}gorin set (see \cite{Brua,Cve,Huang,Me1,C,Me2,LC,Var}). However, it is worth noting here that one did not consider the problem that whether or not there is some proper subset for these sets in which each eigenvalue of a matrix is not included, until Melman in \cite{Me} gave the following Ger\v{s}gorin-type set.

\begin{thm}\label{th1.3} 
Let $A=[a_{ij}]\in \mathbb{C}^{n\times n}$ be a complex matrix. Then
\[\sigma(A) \subseteq \Omega(A) =\bigcup \limits_{i=1}^n \Omega_i(A),\]
where \[\Omega_i(A) = \Gamma_i(A) \backslash \Delta_{i} (A), ~ ~ ~\Delta_{i} (A) =   \bigcup \limits_{j \neq i,\atop j=1}^n \Delta_{ij} (A),\]
and
\[\Delta_{ij} (A)= \left \{z\in \mathbb{C}:|z-a_{jj}| < 2|a_{ji}|-r_j(A) \right \}.\]
Furthermore, $\Omega(A) \subseteq \Gamma(A)$.
\end{thm}

Remark here that there is a typographical error in a
single definition, namely, "$\geq$" instead of "$<$" in the definition of $\Delta_{ij} (A)$ (Theorem 2 in [12]).

Inspired by A. Melman, we in this paper give two new Brauer-type eigenvalue inclusion sets  by considering all but the largest modulus component of an eigenvector and its corresponding
characteristic polynomial equation of a matrix, and by considering  the largest modulus component and the second-largest modulus component. And it is proved that this new Brauer-type sets are better than the Brauer set.

\section{Exclusion sets for the Brauer set}
 In this section,  we present two new Brauer-type sets by excluding two kinds of  Brauer-type exclusion sets from the Brauer set, and the relations between them and the Brauer set are also given.

\begin{thm}\label{th3.1}
Let $A=[a_{ij}]\in \mathbb{C}^{n\times n}$ be a complex matrix. Then
\[\sigma(A)\subseteq \Phi(A)= \bigcup\limits_{j\neq i,\atop i, j=1}^{n} \Phi_{ij}(A),\]
where \[\Phi_{ij}(A)  = \mathcal{K}_{ij}(A) \backslash\mathcal{L}_{i}(A), ~ ~ ~ \mathcal{L}_{i}(A)=\bigcup\limits_{s\neq i, \atop s=1}^{n}\mathcal{L}_{si}(A),\]
\begin{equation} \label{eq8} \mathcal{L}_{si}(A)=\left \{z\in \mathbb{C}:|z -a_{ss}|(|z -a_{ii}|+r_i^s(A)) < (|a_{si}|-r_{s}^{i}(A))|a_{is}| \right \},\end{equation}
and
\[r_t^k(A)=r_t(A)-|a_{tk}|, ~ ~ ~ \forall k\neq t.\]
Furthermore, $\Phi(A)\subseteq\mathcal{K}(A)$.
\end{thm}

\begin{proof} Suppose that $\lambda$ is an eigenvalue of $A$ with a corresponding eigenvector $x=(x_1,x_2,\ldots,x_n)^T$, then
\begin{equation} \label{cha}
Ax =\lambda x
\end{equation} 
holds. Let \[ |x_p| \geq |x_t|\geq\max\limits_{k\neq p, t \atop 1\leq k\leq n} |x_k|,\]
hence $|x_p|>0$. By the proof of the well-known Brauer theorem in \cite{Br}, also see Theorem 2.2 in \cite{Var}, we can easily get
that the $p$-th equality of (\ref{cha}): \begin{eqnarray}\label{meq1} (\lambda -a_{pp})x_p= \sum\limits_{k=1,\atop k\neq p}^{n} a_{pk}x_k\end{eqnarray} gives
\begin{equation}\label{eq10}|\lambda -a_{pp}||x_p|\leq \sum\limits_{k=1,\atop k\neq p}^{n} |a_{pk}||x_k| \leq \sum\limits_{k=1,\atop k\neq p}^{n} |a_{pk}||x_t|= r_p(A)|x_t|,\end{equation}
 and the $t$-th equality of (\ref{cha}):
\begin{eqnarray}\label{meq2}(\lambda -a_{tt})x_t= \sum\limits_{k=1,\atop k\neq t}^{n} a_{tk} x_k\end{eqnarray} gives
\begin{equation}\label{eq12}|\lambda -a_{tt}||x_t|\leq \sum\limits_{k=1,\atop k\neq t}^{n} |a_{tk}||x_k| \leq \sum\limits_{k=1,\atop k\neq t}^{n} |a_{tk}||x_p|= r_t(A)|x_p|.\end{equation}
If $|x_t|=0$, then from (\ref{eq10}), we have $\lambda=a_{pp}$, which implies $\lambda\in \mathcal{K}_{pt}(A)$.
If $|x_t|>0$, by (\ref{eq10}) and (\ref{eq12}), we have
\[|\lambda -a_{pp}||\lambda -a_{tt}|\leq r_p(A)r_t(A),\]
that is, \begin{equation}\label{eq13} \lambda\in \mathcal{K}_{pt}(A).  \end{equation}

On the other hand, for any $s\neq p$, and by the $s$-th equality of (\ref{cha}), we have
\begin{equation}\label{meq6} (\lambda -a_{ss})x_s= \sum\limits_{k=1,\atop k\neq s, p}^{n} a_{sk}x_k+a_{sp}x_p, \end{equation}
then
\begin{equation}\label{eq15}a_{sp}x_p=(\lambda -a_{ss})x_s - \sum\limits_{ k=1,\atop k\neq s, p}^{n} a_{sk}x_k.\end{equation}
Taking absolute values on both sides of (\ref{eq15}) and using the triangle inequality gives
\begin{eqnarray*}
|a_{sp}||x_p|&=&|(\lambda -a_{ss})x_s - \sum\limits_{ k=1,\atop k\neq s, p}^{n} a_{sk}x_k|\\
&\leq&|\lambda -a_{ss}||x_s|+\sum\limits_{k=1,\atop k\neq s, p}^{n} |a_{sk}||x_p|\\
&=&|\lambda -a_{ss}||x_s|+r_{s}^{p}(A)|x_p|,
\end{eqnarray*}
then
\begin{equation}\label{eq16}
(|a_{sp}|-r_{s}^{p}(A))|x_p|\leq |\lambda -a_{ss}||x_s|.
\end{equation}
By the $p$-th equation of (\ref{cha}), we have
\begin{equation}\label{meq5} (\lambda -a_{pp})x_p= \sum\limits_{k=1,\atop k\neq p, s}^{n} a_{pk}x_k+a_{ps}x_s, \end{equation}
then
\begin{equation}\label{eq18}a_{ps}x_s=(\lambda -a_{pp})x_p-\sum\limits_{k=1,\atop k\neq p, s}^{n} a_{pk}x_k.\end{equation}
Taking absolute values on both sides of (\ref{eq18})and using the triangle inequality yields
\begin{eqnarray*}
|a_{ps}||x_s|&=&|(\lambda -a_{pp})x_p-\sum\limits_{k=1,\atop k\neq p, s}^{n} a_{pk}x_k|\\
&\leq&|\lambda -a_{pp}||x_p|+\sum\limits_{k=1,\atop k\neq p, s}^{n} |a_{pk}||x_p|\\
&=&|\lambda -a_{pp}||x_p|+r_{p}^{s}(A)|x_p|,
\end{eqnarray*}
hence
\begin{equation}\label{eq19}
|a_{ps}||x_s|\leq \left(|\lambda -a_{pp}|+r_{p}^{s}(A)\right)|x_p|.
\end{equation}
If $|x_s|>0$, then by (\ref{eq16}) and (\ref{eq19}) gives
\begin{equation}\label{eq20}
|\lambda -a_{ss}|\left(|\lambda -a_{pp}|+r_{p}^{s}(A)\right)\geq(|a_{sp}|-r_{s}^{p}(A))|a_{ps}|,
\end{equation}
that is \begin{equation}\label{eq21}\lambda\notin\mathcal{L}_{sp}(A).\end{equation} Notice that (\ref{eq21}) holds for any $s\neq p$, then \begin{equation}\label{eq22}\lambda\notin\left(\bigcup\limits_{s\neq p,\atop s=1}^{n}\mathcal{L}_{sp}(A)\right)=\mathcal{L}_{p}(A).\end{equation}
From (\ref{eq13}) and (\ref{eq22}), we have
\begin{equation}\label{eq23}
\lambda\in\left(\mathcal{K}_{pt}(A)\backslash\mathcal{L}_{p}(A)\right)=\Phi_{pt}(A).
\end{equation}
Since we do not know which $p$ and $t$ are appropriate to each eigenvalue $\lambda$, we can only conclude that
\begin{equation}\label{eq24}\lambda\in\left(\bigcup\limits_{t\neq p,\atop p,t=1}^{n}\Phi_{pt}(A)\right)=\Phi(A).\end{equation}
On the other hand, if $|x_s|=0$, then from (\ref{eq16}), we have $|a_{sp}|-r_{s}^{p}(A)\leq0$, which implies (\ref{eq20}) holds, and then (\ref{eq24}) holds. Hence \[\sigma(A)\subseteq \Phi(A).\]

In addition, since \[\left(\mathcal{K}_{pt}(A)\backslash\mathcal{L}_{p}(A)\right)\subseteq\mathcal{K}_{pt}(A),\] then \[\Phi(A)\subseteq\mathcal{K}(A).\]
The proof is completed. \qquad
\end{proof}

\begin{rmk}  
 (I)  Theorem \ref{th3.1} shows that for each eigenvalue $\lambda $ of $A$, $\lambda \notin \mathcal{L}_{i}(A)$ for $i\in N$. Note that
$\mathcal{L}_{i}(A)$ is generated by the union of $n-1$ Cassini ovals determined by the elements of $A$, hence $\mathcal{L}_{i}(A)$ is called a Brauer-type exclusion set corresponding to the $(i,j)$-th Brauer Cassini oval $\mathcal{K}_{ij}(A)$.

(II) the computation for $\Phi(A)$ needs $\frac{3n(n-1)}{2}$ Cassini ovals, while the computation for $\mathcal{K}(A)$ needs $\frac{n(n-1)}{2}$ Cassini ovals.
\end{rmk}

{\sc Example\ 3.1} Consider the matrix\[A= \left[
\begin{array}{cccc}
  14   &0.01\textbf{i}     &0     &18-2\textbf{i} \\
  0    &9    &4+\textbf{i}   &0  \\
  0.01+\textbf{i}  &2+\textbf{i}   &11    &0\\
  19+\textbf{i} & 0    &0.1+\textbf{i}   &10
\end{array} \right].\]
The sets $\Phi(A)$ in Theorem \ref{th3.1} is drawn in
Figure 1. And the exact eigenvalues of $A$ are plotted with asterisks. From Figure 1, we conclude that $\Phi(A)$ locate the eigenvalues of $A$ more precisely than $\mathcal{K}(A)$.


Note that Theorem \ref{th3.1} is obtained by considering all but the largest modulus component of an eigenvector and its corresponding characteristic polynomial
 equation, which needs much computations. To reduce its computations, we next give another Brauer-type set by considering only the largest modulus component and the second-largest modulus component of an eigenvector and its corresponding characteristic polynomial equation.

\begin{thm}\label{th3.2}
Let $A=[a_{ij}]\in \mathbb{C}^{n\times n}$ be a complex matrix. Then
\[\sigma(A)\subseteq \Theta(A)= \bigcup\limits_{j\neq i,\atop i,j=1}^{n} \Theta_{ij}(A),\]
where \[\Theta_{ij}(A)  = \mathcal{K}_{ij}(A) \backslash \Lambda_{ij}(A) ,\]
and
\begin{equation} \label{eq25} \Lambda_{ij}(A)= \left \{z\in \mathbb{C}:(|\lambda -a_{ii}|+r_i^j(A)) (|\lambda -a_{jj}|+r_j^i(A)) < |a_{ij}||a_{ji}| \right \}.\end{equation}
Furthermore, $\Theta(A)  \subseteq \mathcal{K}(A)$.
\end{thm}

\begin{proof} Suppose that $\lambda$ is an eigenvalue of $A$ with a corresponding eigenvector $x=(x_1,x_2,\ldots,x_n)^T$, then
(\ref{cha}) holds. Let \[ |x_p| \geq |x_t|=\max\limits_{k\neq p,t\atop 1\leq k\leq n } |x_k|,\]
then $|x_p|>0$. Similar to the proof of Theorem \ref{th3.1},
\begin{equation}\label{eq26} \lambda\in \mathcal{K}_{pt}(A) \end{equation} can be easily obtained.

On the other hand, (\ref{meq1}) and (\ref{meq2}) can be rewritten respectively as
\begin{equation}\label{eq27}(\lambda -a_{pp})x_p-\sum\limits_{k=1,\atop k\neq p, t}^{n} a_{pk}x_k= a_{pt}x_t\end{equation}
and
\begin{equation}\label{eq28}(\lambda -a_{tt})x_t- \sum\limits_{k=1,\atop k\neq t, p}^{n} a_{tk} x_k=a_{tp}x_p.\end{equation}
Taking absolute values on both sides of (\ref{eq27}) and (\ref{eq28}), and using the triangle inequality yields
\begin{equation}\label{eq29}
|a_{pt}||x_t| \leq |\lambda -a_{pp}||x_p|+\sum\limits_{k=1,\atop k\neq p, t}^{n} |a_{pk}||x_k|
\leq (|\lambda -a_{pp}|+r_p^t(A)) |x_p|\end{equation}
and
\begin{equation}\label{eq30}|a_{tp}||x_p| \leq |\lambda -a_{tt}||x_t|+ \sum\limits_{k=1,\atop k\neq t, p}^{n} |a_{tk}| |x_k|\leq (|\lambda -a_{tt}|+r_t^p(A)) |x_t|.\end{equation}
If $|x_t|\neq 0 $, then multiplying (\ref{eq29}) and (\ref{eq30}) gives
\[ |a_{pt}||a_{tp}||x_t||x_p|\leq (|\lambda -a_{pp}|+r_p^t(A)) (|\lambda -a_{tt}|+r_t^p(A)) |x_p||x_t|.\]
that is
\begin{equation}\label{eq31}  (|\lambda -a_{pp}|+r_p^t(A)) (|\lambda -a_{tt}|+r_t^p(A))\geq |a_{pt}||a_{tp}|,\end{equation}
which implies that
\begin{equation}\label{eq32}\lambda\notin \Lambda_{pt}(A). \end{equation}
If $|x_t|= 0 $, then by (\ref{eq30}) we have $|a_{tp}| =0 $, which also leads to $\lambda\notin \Lambda_{pt}(A)$. Furthermore,
from (\ref{eq26}) and (\ref{eq32}), we have
\begin{equation}\label{eq33}\lambda\in \left(\mathcal{K}_{pt}(A) \backslash \Lambda_{pt}(A)\right) = \Theta_{pt}(A).\end{equation}
Since we do not know which $p$ and $t$ are corresponding to each eigenvalue $\lambda$, then we can only get that \[  \lambda\in  \left(\bigcup\limits_{t\neq p} \Theta_{pt}(A) \right)= \Theta(A).\]
Hence \[ \sigma(A)\subseteq \Theta(A).\]

In addition, since \[\left(\mathcal{K}_{pt}(A)\backslash \Lambda_{pt}(A)\right)\subseteq\mathcal{K}_{pt}(A),\] then \[\Phi(A)\subseteq \mathcal{K}(A)\] can be easily obtained. The conclusion follows.
\end{proof}

\begin{rmk}
 (I) Note that the computation for $\Theta(A)$ needs $n(n-1)$ Cassini ovals, which is obviously less than that of $\Phi(A)$.

(II) By lots of numerical examples we find that $\Phi(A)\subset\Theta(A)$ in most cases, and the worst case is $\Phi(A)=\Theta(A)$. We here give a conjecture that
\[ \Phi(A) \subseteq \Theta(A).\]
\end{rmk}

{\sc Example\ 3.2} Consider again the matrix $A$ in Example 3.1. The set $\Theta(A)$ in Theorem \ref{th3.2} is drawn in  Figure 2, and the exact eigenvalues of $A$ are plotted with asterisks. From Figure 2, it is not difficult to see that $\Theta(A)$ can also locate the eigenvalues of $A$ more precisely than $\mathcal{K}(A)$, but comparing   Figure 1 with  Figure 2, we find that $\Phi(A)\subset\Theta(A)$.


Just as the Ger\v{s}gorin disk theorem leads to
 the condition of strict diagonal dominance \cite{Lji}, we next give two sufficient criterions for the non-singularity of complex matrices by Theorem \ref{th3.1} and Theorem \ref{th3.2}.

\begin{corollary}\label{cor3.1}
Let $A=[a_{ij}]\in \mathbb{C}^{n\times n}$ be a complex matrix. Then $A$ is non-singular if for each $i\in\{1,\cdots,n\}$, $j\in\{1,\cdots,n\}$, and $j\neq i$, either
\[|a_{ii}||a_{jj}|> r_i(A)r_j(A)\]
or
\[|a_{ss}|(|a_{ii}|+r_i^s(A)) < (|a_{si}|-r_{s}^{i}(A))|a_{is}|\]for some $s\neq i$ and $s\in\{1,\cdots,n\}$.
\end{corollary}

\begin{corollary}\label{cor3.2}
Let $A=[a_{ij}]\in \mathbb{C}^{n\times n}$ be a complex matrix. Then $A$ is non-singular if for each $i\in\{1,\cdots,n\}$, $j\in\{1,\cdots,n\}$, and $j\neq i$, either
\[|a_{ii}||a_{jj}|> r_i(A)r_j(A)\]
or
\[(|a_{ii}|+r_i^j(A)) (|a_{jj}|+r_j^i(A)) < |a_{ij}||a_{ji}|.\]
\end{corollary}

\section{Conclusion}
In this paper, two new Brauer-type sets $\Phi(A)$ and $\Theta(A)$ are given by excluding its corresponding Brauer-type
exclusion sets, respectively. To investigate the relations between this two new Brauer-type sets and the Brauer set, we compare them with each other and obtain a novel result. Actually, by the similar method,
there are many eigenvalue inclusion sets, such as, the sets in \cite{Cve,Huang,Me1,Me2,LC}, from which we can exclude their corresponding exclusion subsets to provide more precise eigenvalue inclusion sets.

\section*{Acknowledgements}
This work is partly supported by National Natural Science
Foundations of China (11601473) and CAS
"Light of West China" Program.







\end{document}